\documentclass[12pt]{amsart}

\usepackage{amsthm,amssymb,amsmath}
\usepackage{rotating}
\usepackage[cm]{fullpage}
\usepackage{url}

\theoremstyle{plain}
\newtheorem{theorem}{Theorem}
\newtheorem{lemma}{Lemma}
\newtheorem{corollary}{Corollary}
\newtheorem{conjecture}{Conjecture}
\newtheorem{problem}{Problem}
\theoremstyle{remark}
\newtheorem{example}{Example}

\renewcommand{\leq}{\leqslant}
\renewcommand{\ge}{\geqslant}

\def\eref#1{$(\ref{#1})$}
\def\egref#1{Example~$\ref{#1}$}
\def\sref#1{\S$\ref{#1}$}
\def\lref#1{Lemma~$\ref{#1}$}
\def\tref#1{Theorem~$\ref{#1}$}
\def\cref#1{Corollary~$\ref{#1}$}
\def\cjref#1{Conjecture~$\ref{#1}$}
\def\Tref#1{Table~$\ref{#1}$}

\def\pref#1{Problem~$\ref{#1}$}

\def\id{{\ensuremath\varepsilon}}

\newcommand{\ldiv}{\backslash}
\newcommand{\rdiv}{/}

\title{Loops with exponent three in all isotopes}

\author{Michael Kinyon}
\address{Department of Mathematics \\
2280 S. Vine St. \\
University of Denver \\ Denver, CO 80208 USA}
\email{\url{michael.kinyon@du.edu}}

\author{Ian M.\ Wanless$^{\dagger}$}
\address{School of Mathematical Sciences\\
Monash University\\
Vic 3800, Australia}
\email{\url{ian.wanless@monash.edu}}

\thanks{${}^{\dagger}$Research supported by ARC grant DP0662946.}

\date{\today}

\begin{document}

\begin{abstract}
It was shown by van Rees \cite{vR} that a latin square of order $n$
has at most $n^2(n-1)/18$ latin subsquares of order $3$.  He
conjectured that this bound is only achieved if $n$ is a power of
$3$. We show that it can only be achieved if $n\equiv3\bmod6$.  We also
state several conditions that are equivalent to achieving the van Rees
bound. One of these is that the Cayley table of a loop achieves the
van Rees bound if and only if every loop isotope has exponent $3$. We
call such loops \emph{van Rees loops} and show that they form an
equationally defined variety.

We also show that (1) In a van Rees loop, any subloop of index 3 is
normal, (2) There are exactly 6 nonassociative van Rees loops of order
$27$ with a non-trivial nucleus and at least 1 with all nuclei trivial,
(3) Every commutative van Rees loop has the weak inverse property and
(4) For each van Rees loop there is an associated family of Steiner 
quasigroups.
\end{abstract}

\maketitle

\section{Introduction}

The background prerequisites for this paper are the basic theories
of quasigroups and loops \cite{Bruck,Pf} and of latin squares \cite{DK}.
Our results and proofs are a mix of the algebraic (quasigroups and loops)
and combinatorial (latin squares) perspectives.

Throughout this paper $Q=(Q,\cdot)$ will denote a quasigroup or, in case
there is a neutral element $\id$, a loop. We denote by $L$ a latin square
obtained from the (unbordered) Cayley table of $Q$ (once an arbitrary ordering
of the elements of $Q$ has been fixed). For $x\in Q$, we define the \emph{left}
and \emph{right translations} $L_x:Q\rightarrow Q$ and $R_x:Q\rightarrow Q$ by,
respectively, $L_x(y) = xy$ and $R_x(y) = yx$ for $y\in Q$.

As usual, the latin square properties in which we are interested are
those invariant, at the very least, under \emph{isotopy}, that is,
under permutations of the rows, of the columns and of the
symbols. From the quasigroup perspective, we are interested in
properties that hold in all \emph{loop isotopes}. Given a quasigroup
or loop $(Q,\cdot)$ and fixed elements $a,b\in Q$, we can define a
loop $Q_{a,b} = (Q,\circ)$ with neutral element $\id = ba$ by $x\circ
y = R_a^{-1}(x)\cdot L_b^{-1}(y)$. Further, all loops isotopic to (the
quasigroup or loop) $(Q,\cdot)$ are isomorphic to isotopes of this
form \cite{Bruck}.  We will denote left and right translations in the
loop isotope $Q_{a,b} = (Q,\circ)$ by $L_x^{\circ}$ and $R_x^{\circ}$,
respectively.

Since the notion of power of an element is not unambiguously defined in
loops, in general it does not make sense to speak of the exponent of a loop.
However, the loops we consider in this paper will satisfy the identity
$xx\cdot x = x\cdot xx = \id$, and in this case we say that the loop has
\emph{exponent} $3$.

A {\em subsquare} of a latin square is a submatrix that is itself a
latin square.  In \cite{vR}, van Rees showed that no latin square of
order $n$ has more than $n^2(n-1)/18$ subsquares of order $3$.

Our main result is the following.

\begin{theorem}\label{t:main}
Suppose $(Q,\cdot)$ is a quasigroup of order $n$ with associated
latin square $L$. The following conditions are equivalent:
\begin{enumerate}
\item \label{i:bnd}
$L$ has $n^2(n-1)/18$ subsquares of order $3$.
\item \label{i:eachpair}
For any two occurrences of the same symbol in $L$, there is
a subsquare of order $3$ containing those two occurrences.
\item \label{i:eachcell}
Every cell in $L$ is in $(n-1)/2$ subsquares of order $3$.
\item \label{i:exp3}
Every loop isotopic to $Q$ has exponent $3$.
\item \label{i:2transl}
For any distinct $x,y\in Q$, $L_x^{-1}L_y$ and $R_x^{-1}R_y$
are regular permutations of order $3$.
\item \label{i:translate}
In any loop isotope $(Q,\circ)$ and for each nonidentity element $x\in Q$,
$L_x^{\circ}$ and $R_x^{\circ}$ are regular permutations of order $3$.
\item \label{i:subloop}
In any loop isotopic to $Q$, there are $(n-1)/2$ subloops
of order $3$.
\end{enumerate}
\end{theorem}

Note that (\ref{i:exp3}), (\ref{i:translate}) and (\ref{i:subloop})
are loop isotope conditions, (\ref{i:2transl}) is a quasigroup (or loop)
condition and (\ref{i:bnd}), (\ref{i:eachpair}) and
(\ref{i:eachcell}) are latin square conditions.

For reasons that will become apparent later,
we call any quasigroup or loop satisfying \tref{t:main}(\ref{i:exp3}) a
\emph{van Rees quasigroup} or \emph{van Rees loop}, respectively,
and any latin square satisfying \tref{t:main}(\ref{i:bnd}) a
\emph{van Rees latin square}.

This paper concentrates on the number of subsquares of order $3$.  For
bounds on the number of subsquares of other orders, see
\cite{BCW,BSW,BVW}. The situation with order $2$ subsquares is fairly well
understood. A latin square of order $n$ can have at most
$n^2(n-1)/4$ subsquares of order $2$, with equality being
achieved precisely when the latin square is isotopic to an elementary
abelian $2$-group. This result, and equivalent statements, has been
rediscovered a number of times (e.g. \cite{Cam72,HW,TD}).
The property of having $n^2(n-1)/4$ subsquares of order $2$
is equivalent to all loop isotopes having exponent $2$. The subject
of this paper is arguably the first non-trivial case regarding
upper bounds for number of subsquares.
There are also lower bounds \cite{Bar13,BCW} on the maximum
number of subsquares of
order $2$ in cases when the upper bound is not achieved.
We do not investigate analogous lower bounds for order $3$ subsquares
here, although it would make a worthy subject for future investigations.


The outline of the paper is as follows. In the next section, after
some preliminaries, we give the proof of \tref{t:main},
followed by some immediate consequences.  In particular, quasigroups
and loops satisfying the conditions of the theorem turn out to form
varieties, not only of quasigroups (\tref{t:variety}), but also
of magmas, that is, sets with a single binary operation 
(\tref{t:variety2}).  In \sref{s:nonexamples} we give examples showing
that the various conditions defining van Rees loops are, in fact,
necessary. In \sref{s:vR}, we examine a conjecture of van Rees
regarding the possible orders of van Rees latin squares, and show that
such a square has order congruent to $3$ (mod $6$) (\tref{t:3mod6}).  
In \sref{s:constructions}, we look at several
examples of nonassociative van Rees loops.  In \sref{s:steiner}, we
show that on the underlying set of a van Rees loop, there is a natural
Steiner quasigroup structure (\tref{t:steiner}).

\section{Proof of the Main Theorem}\label{s:proof}

We will need the following two elementary results, which are well known
and easy to prove.

\begin{lemma}\label{l:basic}
Let $L$ be a latin square.
\begin{enumerate}
\item \label{i:intersect}
If two subsquares of a latin square $L$ have nontrivial intersection,
that intersection is itself a subsquare.
\item \label{i:half}
If $S$ is a subsquare of a latin square $L$ and $S\ne L$ then the
order of $S$ cannot exceed half the order of $L$.
\end{enumerate}
\end{lemma}

\begin{lemma}\label{l:subloop}
Suppose $Q$ is a loop with associated latin square $L$ and neutral
element $\id$.  If $S$ is a subsquare of $L$ including the cell
$(\id,\id)$ then the set of symbols occurring in $S$ is a subloop of $Q$.
\end{lemma}

\begin{proof}[Proof of Theorem \ref{t:main}]
Let $S_1, S_2$ be subsquares of order $3$ containing two different
occurrences $u,v$ of the same symbol $x$ in $L$. Since $S_1\cap S_2$
has order at least $2$, \lref{l:basic} implies that $S_1 = S_2$. Thus
any two occurrences of a symbol are in \emph{at most} one subsquare
of order $3$.

Now, there are $n$ choices for the symbol $x$ and $\binom{n}{2}$
choices for the occurrences $u,v$. Each subsquare of order $3$ contains
$3{\binom{3}{2}} = 9$ pairs of entries with the same symbol. It follows that
$L$ can have at most
\[
\frac{n{\binom{n}{2}}}{3{\binom{3}{2}}} = \dfrac{1}{18} n^2(n-1)
\]
subsquares of order $3$, with this bound being achieved if and only if
condition (\ref{i:eachpair}) holds. That is,
$(\ref{i:bnd})\Longleftrightarrow(\ref{i:eachpair})$

(\ref{i:eachpair})$\Longleftrightarrow$(\ref{i:eachcell}):
For a fixed choice of an occurrence $u$ of a symbol, we have $n-1$ choices
for a different occurrence $v$, and each
subsquare of order $3$ that includes $u$ also includes two options
for $v$.

(\ref{i:exp3})$\implies$(\ref{i:2transl}):
If the loop isotope $Q_{a,b}$ has exponent $3$, then for all $x\in Q$,
\[
ba = xa\circ (xa\circ xa)
= R_a^{-1}(xa)\cdot L_b^{-1}(R_a^{-1}(xa) \cdot L_b^{-1}(xa))
= L_x L_b^{-1} L_x L_b^{-1} L_x (a),
\]
and similarly, $ba = (bx\circ bx)\circ bx = R_x R_a^{-1} R_x R_a^{-1} R_x(b)$.
Thus $(L_b^{-1} L_x)^3(a) = a$ and $(R_a^{-1} R_x)^3(b) = b$. We conclude
that if every loop isotope of $Q$ has exponent $3$, then each
$L_b^{-1} L_x$ and each $R_a^{-1} R_x$ is a regular permutation of order $3$.

(\ref{i:2transl})$\Longleftrightarrow$(\ref{i:translate}):
In the loop isotope $Q_{a,b}$, left translations are given by
$L_x^{\circ} = L_{R_a^{-1}(x)} L_b^{-1}$ and right translations are
given by $R_x^{\circ} = R_{L_b^{-1}(x)} R_a^{-1}$. Thus
$L_x^{\circ}$ and $R_x^{\circ}$ are each regular of order $3$
for \emph{all} $x\in Q$ if and only if $L_x L_b^{-1}$
and $R_x R_a^{-1}$ are each regular of order $3$ for \emph{all}
$x\in Q$. Universally quantifying $a$ and $b$, we have the desired
equivalence.

(\ref{i:2transl})$\implies$(\ref{i:eachpair}):
Suppose $s$ is a symbol in $L$ occurring in the distinct cells
$(a,b)$ and $(c,d)$ so that $s = ab = cd$ where $a\neq c$.
Set $t = ad$ and $u = cb$. Since $L_a L_c^{-1}$ is regular
of order $3$, $t = L_a L_c^{-1}(s) = L_a L_c^{-1} L_a L_c^{-1} (cb)
= L_c L_a^{-1}(u) = c\cdot L_a^{-1}(u)$.
Similarly, since $R_d R_b^{-1}$ is regular of order $3$,
$t = R_d R_b^{-1}(s) = R_d R_b^{-1} R_d R_b^{-1} (cb)
= R_b R_d^{-1}(u) = R_d^{-1}(u)\cdot b$. Now set
$e = R_d^{-1}(u) = R_b^{-1}(t)$ and $f = L_a^{-1}(u) = L_c^{-1}(t)$.
Note that $L_e^{-1}(u) = d$ and $L_e^{-1}(t) = b$.
Since $L_c L_e^{-1}$ is regular of order $3$, we have
$s = cd = L_c L_e^{-1}(u) = L_c L_e^{-1} L_c(b)
= L_c L_e^{-1} L_c L_e^{-1}(t) = L_e L_c^{-1}(t) = ef$.
From these calculations, we deduce that $L$ has a subsquare of order
$3$ containing the two occurrences of $S$, namely the subsquare labeled
by the rows $a,c,e$ and the columns $b,d,f$.

(\ref{i:eachpair})$\implies$(\ref{i:exp3}):
Consider the loop isotope
$Q_{a,b}$ with neutral element $ba$.  Let $s\in
Q\setminus\{ba\}$. By assumption, cells $(ba,s)$ and
$(s,ba)$ are in a subsquare $S$ of order $3$, and hence by
\lref{l:basic} are not in a subsquare of order $2$. Thus $s^2\ne ba$
and the symbols in $S$ must be $ba,s,s^2$. Moreover, these symbols
form a subloop of $Q_{a,b}$ by \lref{l:subloop}. As $s$ was an
arbitrary non-identity element, we conclude that $Q_{a,b}$ has exponent
$3$, that is, condition (\ref{i:exp3}) holds.

(\ref{i:subloop})$\implies$(\ref{i:eachcell}):
Let $(a,b)$ be any cell in $L$.
In the loop isotope $Q_{b,a}$ there are, by assumption,
exactly $(n-1)/2$ subloops of order $3$, each
of which corresponds to a different subsquare of order $3$ including the
cell $(a,b)$.

(\ref{i:eachcell})$\implies$(\ref{i:subloop}):
In all loop isotopes $Q_{a,b}$, the
cell $(ba,ba)$ is in exactly $(n-1)/2$ subsquares of order $3$, each
of which corresponds to a different subloop of order $3$, by \lref{l:subloop}.
\end{proof}

\begin{corollary}\label{c:parastrophe}
Any parastrophe of a van Rees quasigroup is also a van Rees quasigroup.
\end{corollary}

\begin{proof}
Parastrophy preserves the number of subsquares of order $3$.
\end{proof}

The permutations $L_x^{-1}L_y$ and $R_x^{-1}R_y$ that appear in
\tref{t:main}(\ref{i:2transl}) are important in the combinatorics of
latin squares.  Their expected structure in a randomly chosen latin
square was studied in \cite{CGW}. Meanwhile
\cite{BMW,MW,perf,cyclatom} examine the case when $L_x^{-1}L_y$ and
$R_x^{-1}R_y$ consist of a single cycle, regardless of the choice of
distinct $x,y$.

We recall the universal algebraic definition of a
quasigroup $(Q;\cdot,\ldiv,\rdiv)$ is a set
$Q$ together with three operations $\cdot, \ldiv, \rdiv : Q\times Q\to Q$
satisfying the identities $x\ldiv (xy) = x(x\ldiv y) = y$ and
$(xy)\rdiv y = (x\rdiv y)y = x$.

\begin{theorem}\label{t:variety}
The class of van Rees quasigroups \textup{[}loops\textup{]} $(Q,\cdot,\ldiv,\rdiv)$
forms a variety of quasigroups $[$loops$]$ defined by the identities
\begin{align*}
x(y\ldiv (xz)) &= y(x\ldiv (yz))  \tag{vR1}\\
((xy)\rdiv z)y &= ((xz)\rdiv y)z.  \tag{vR2}
\end{align*}
\end{theorem}

\begin{proof}
This is just Theorem \ref{t:main}(\ref{i:2transl}) written explicitly in
terms of the divisions $\ldiv$ and $\rdiv$.
\end{proof}

The algebraic advantage of the three operation definition of quasigroups
over the one operation definition is that homomorphic images under the
latter definition need not be quasigroups. However, in this case, we
can also view van Rees loops as varieties of magmas $(Q,\cdot)$.

\begin{theorem}\label{t:variety2}
The class of van Rees loops $(Q,\cdot)$ forms a variety of magmas
(with neutral element) defined by the identities
\begin{align*}
x(x\cdot xy) &= y \tag{vRL1} \\
(xy\cdot y)y &= x \tag{vRL2} \\
x\cdot y(y\cdot xz) &= y\cdot x(x\cdot yz) \tag{vRL3} \\
(xy\cdot z)z\cdot y &= (xz\cdot y)y\cdot z.\tag{vRL4}
\end{align*}
\end{theorem}

\begin{proof}
If $(Q,\cdot,\ldiv,\rdiv)$ satisfies (vR1), then taking $x = \id$,
we obtain (vRL1), and similarly, (vR2) implies (vRL2). Now
$x\ldiv y = x\cdot xy$ and $x\rdiv y = xy\cdot y$, and so (vRL3) and (vRL4)
are just (vR1) and (vR2) rewritten. Conversely, assume $(Q,\cdot,\ldiv,/)$
satisfies (vRL1)--(vRL4). Define $x\ldiv y = x(xy)$, and observe that (vRL1) implies
$x\ldiv (xy) = x(x\ldiv y) = y$. Similarly defining
$x\rdiv y = (xy)y$, (vRL2) implies $(xy)\rdiv y = (x\rdiv y)y = x$. Thus
$(Q,\cdot,\ldiv,\rdiv )$ is a quasigroup, and then (vR1) and (vR2) are
just (vRL3) and (vRL4), respectively, rewritten.
\end{proof}

The following is immediate from either Theorem \ref{t:variety} or \ref{t:variety2}.

\begin{corollary}\label{c:sublpclosed}
The class of van Rees quasigroups \textup{[}loops\textup{]} is closed
under taking homomorphic images, subquasigroups \textup{[}subloops\textup{]} and direct products. 
Any subsquare of a van Rees latin square is a van Rees latin square.
\end{corollary}

\begin{proof}
The first statement follows from Birkhoff's theorem on the equivalence
of varieties of universal algebras and equational classes. The second
statement follows because for any subsquare, there is a loop isotope
that turns the subsquare into a subloop, in the sense of \lref{l:subloop}.
\end{proof}

\section{Nonexamples}
\label{s:nonexamples}

We now consider some examples of loops which are not van Rees loops
to show that the various defining conditions are necessary.
These examples were found by a mixture of the finite model builder
\textsc{Mace4} \cite{McCune} and by home-grown software.

\begin{example}\label{ex:notinv}
Having exponent $3$ is not, \emph{a priori}, an isotopically invariant property
of a loop. The smallest counterexample is given in \Tref{t:exp3}.
\begin{table}[htb]
\centering
\small\texttt{
\begin{tabular}{c|ccccccc}
$\cdot$ & \id & a & b & c & d & e & f\\
\hline
    \id & \id & a & b & c & d & e & f \\
    a & a & b & \id & e & f & c & d \\
    b & b & \id & a & f & e & d & c \\
    c & c & e & f & d & \id & a & b \\
    d & d & f & e & \id & c & b & a \\
    e & e & c & d & a & b & f & \id \\
    f & f & d & c & b & a & \id & e
\end{tabular}}
\smallskip
\caption{\label{t:exp3}A loop of exponent $3$.}
\end{table}
This loop has exponent $3$, is commutative and satisfies the identity $x(yx)^2 = y^2$,
which (for loops of exponent $3$) is the weak inverse property. In fact, this loop is
isotopic to the Steiner quasigroup of order $7$. In addition, it turns out that every
isotope of the loop is power-associative. However, to see that it fails to be
a van Rees loop, observe, for instance, that there is a subsquare of order $2$
formed by rows $a,b$ and columns $c,d$. Then apply \cref{c:sublpclosed}.
\end{example}

\begin{example}\label{ex:regnot}
A loop of exponent $3$ can have all left and right translations being regular
permutations of order $3$, but still not have the van Rees property.
Put another way, the identities (vRL3) and (vRL4) in Theorem \ref{t:variety2}
cannot be dispensed with. The smallest example showing this is given in
Table \ref{t:reg}.
\begin{table}[htb]
\centering
\small\texttt{
\begin{tabular}{c|ccccccccc}
$\cdot$ & \id & a & b & c & d & e & f & g & h\\
\hline
    \id & \id & a & b & c & d & e & f & g & h \\
    a & a & b & \id & d & h & f & g & e & c \\
    b & b & \id & a & e & f & g & h & c & d \\
    c & c & g & h & f & a & b & \id & d & e \\
    d & d & h & c & g & e & \id & a & b & f \\
    e & e & c & g & h & \id & d & b & f & a \\
    f & f & d & e & \id & g & h & c & a & b \\
    g & g & e & f & a & b & c & d & h & \id \\
    h & h & f & d & b & c & a & e & \id & g
\end{tabular}}
\smallskip
\caption{\label{t:reg}A loop of exponent $3$ with left and right
translations regular of order $3$.}
\end{table}
To see that this is not a van Rees loop, observe that the occurrences of
$d$ in the cells $(\id,d)$ and $(a,c)$ do not lie in a subsquare of order $3$.
\end{example}

\begin{example}\label{ex:onesided}
Similarly, in Theorem \ref{t:variety}, the conditions (vR1) and (vR2)
are independent. In other words, for a quasigroup or loop to have the
van Rees property, it is not sufficient that, say, $L_x^{-1}L_y$ be regular of order $3$ for
all distinct $x,y\in Q$. An example is the loop defined on $\mathbb{Z}_3\times \mathbb{Z}_3$ by
\[
(x,a)(y,b) = (x+y,a+b+x^2 y).
\]
This is one of the three nonassociative conjugacy closed loops of order $9$, known as
CCLoop(9,3) in the LOOPS package \cite{LOOPS} for GAP \cite{GAP}. This loop does not even have a
well-defined exponent, for although it satisfies the identity $x(xx) = \id$, it does not
satisfy $(xx)x = \id$. Indeed, if $(xx)x=x(xx)$ for all $x$, then the loop would be power-associative,
but the only power-associative conjugacy closed loops of order $9$ are groups \cite{KK}.
\end{example}

Regarding axiomatic considerations, we have not been able to resolve the following.

\begin{problem}\label{p:vR1notvR2}
  Does there exist a loop $(Q,\cdot,\ldiv,\rdiv)$ satisfying the
  identities \emph{(vR1)} and $(xx)x = x(xx) = \id$, but not \emph{(vR2)}?
\end{problem}

\pref{p:vR1notvR2} is motivated by Example~\ref{ex:onesided}. Put
another way, the problem asks for a loop of exponent 3 in which all
permutations $L_x^{-1} L_y$ are regular of order $3$, but in which
some permutation $R_a^{-1} R_b$ does not have order $3$. A stronger
requirement is the following.

\begin{problem}\label{p:vRL123not4}
  Does there exist a loop $(Q,\cdot)$ satisfying the identities
  \emph{(vRL1)}, \emph{(vRL2)} and \emph{(vRL3)}, but not \emph{(vRL4)}?
\end{problem}

Again rephrasing, \pref{p:vRL123not4} asks for a loop (necessarily of
exponent $3$) in which all permutations $L_x^{-1} L_y$ and all
permutations $R_x$ are regular of order $3$, but in which some
permutation $R_a^{-1} R_b$ does not have order $3$.

\section{van Rees' conjecture}
\label{s:vR}

In \cite{vR}, van Rees conjectured that a latin square of order
$n$ cannot have $n^2(n-1)/18$ subsquares of order $3$ unless $n$ is a
power of $3$.  This conjecture provided the original motivation for
the present paper. In our terminology it can be stated like this:

\begin{conjecture}\label{c:vR}
If $L$ is a van Rees latin square, then the order of $L$ is a power of $3$.
\end{conjecture}

\begin{example}\label{ex:notreg2}
Referring to Example \ref{ex:regnot}, it is tempting to make a stronger
conjecture than van Rees', namely that any loop of exponent $3$ in which
every left and right translation is a regular permutation of order $3$ has
order a power of $3$. However, this is not true, as the loop given in
Table \ref{t:notreg2} shows.
\begin{table}[htb]
\centering
\small\texttt{
\begin{tabular}{c|ccccccccccccccc}
$\cdot$ & \id& a& b& c& d& e& f& g& h&i&j&k&l&m&n\\
\hline
\id & \id&a&b&c&d&e&f&g&h&i&j&k&l&m&n\\
a & a&b&\id&e&g&f&c&h&d&n&k&m&i&j&l\\
b & b&\id&a&f&h&c&e&d&g&l&m&j&n&k&i\\
c & c&i&j&d&\id&n&m&f&e&k&l&a&b&g&h\\
d & d&k&l&\id&c&h&g&m&n&a&b&i&j&f&e\\
e & e&m&h&g&a&l&j&i&k&c&n&b&\id&d&f\\
f & f&g&n&h&b&i&k&l&j&m&c&\id&a&e&d\\
g & g&n&f&a&e&k&i&j&l&b&\id&d&m&h&c\\
h & h&e&m&b&f&j&l&k&i&\id&a&n&d&c&g\\
i & i&j&c&m&l&a&d&n&\id&h&e&g&f&b&k\\
j & j&c&i&n&k&d&b&\id&m&f&g&e&h&l&a\\
k & k&l&d&j&n&m&\id&a&c&e&h&f&g&i&b\\
l & l&d&k&i&m&\id&n&c&b&g&f&h&e&a&j\\
m & m&h&e&l&i&g&a&b&f&j&d&c&k&n&\id\\
n & n&f&g&k&j&b&h&e&a&d&i&l&c&\id&m
\end{tabular}}
\smallskip
\caption{\label{t:notreg2}Another loop of exponent $3$ with left and right
translations regular of order $3$.}
\end{table}
If this were a van Rees loop, it would have $175$ subsquares
of order $3$, but it has only $24$ such subsquares. However,
it does satisfy $L_x^3 = R_x^3 = 1$ for all $x$. In addition,
this particular loop has the inverse property. 
\end{example}

\begin{lemma}\label{l:exp3}
A finite loop of exponent $3$ has odd order.
\end{lemma}

\begin{proof}
Let $Q$ be a loop of exponent $3$, and consider the mapping
$Q\to Q; x\mapsto x^2 = x^{-1}$. This map is an involution of
the set $Q^* = Q\ldiv \{1\}$ of nonidentity elements of $Q$, and
does not fix any elements of $Q^*$. Thus $|Q^*|$ is even, and hence
$|Q|$ is odd.
\end{proof}

We have been unable to settle Conjecture \ref{c:vR}, but we can show:

\begin{theorem}\label{t:3mod6}
A van Rees quasigroup or van Rees latin square has order $n\equiv3\bmod6$.
\end{theorem}

\begin{proof}
If there exists a van Rees quasigroup of order $n$, then $n$ is odd by 
\lref{l:exp3} and $n$ is divisible by $3$ by \tref{t:main}(\ref{i:translate}).
\end{proof}

By computer search, we have established that there are no van Rees
loops of orders $15$ or $21$.  Hence the smallest
possible order for a counterexample to van Rees' conjecture is $33$.
Also, given \cref{c:sublpclosed}, we conclude that every subloop of a
van Rees loop has order $m\equiv 3\mod 6$ where $m$ is a power of $3$ or
$m\ge 33$. Regarding the order of subloops, we can also say this:

\begin{theorem}\label{t:nobigsubsq} $\mbox{ }$
\begin{enumerate}
\item If $S$ is a proper subsquare of a van Rees latin square $L$, then the
  order of $S$ is no more than one third of the order of $L$;
\item If $P$ is a proper subquasigroup \textup{[}subloop\textup{]} 
of a van Rees quasigroup \textup{[}loop\textup{]} $Q$,
then $|P| \leq \frac{1}{3} |Q|$.
\end{enumerate}
\end{theorem}

\begin{proof}
(1) Since $S$ is a proper subsquare we can find a row of $L$ indexed by, say,
$x$ that intersects $S$ and another row indexed by, say, $y$ that does not.
The permutation $L_x^{-1}L_y$ consists of cycles of length $3$. Consider such
a cycle $C$ containing a symbol $s$ from inside $S$. Neither
the image nor preimage of $s$ in $C$ can be symbols from $S$.
Therefore there are at least as many $3$-cycles in $L_x^{-1}L_y$ as
there are symbols in $S$. Thus (1) follows, and then (2) follows from (1).
\end{proof}

In the case when the bound in \tref{t:nobigsubsq} is achieved, we can say more.

\begin{theorem}\label{t:normal}
If $S$ is a subloop of index $3$ in a van Rees loop $Q$,
then $S$ is normal in $Q$.
\end{theorem}

\begin{proof}
For $q\in Q\setminus S$ define
$T_q=\{x\ldiv q:x\in S\}$ and $U_q=\{q\rdiv x:x\in S\}$.
Let $\Sigma_0$ be the subsquare $\{(s_1,s_2,s_1\circ s_2):s_1,s_2\in S\}$.


Suppose $r_1,c_1\in S$. By \tref{t:main} there is a subsquare $\Sigma_3$
of order $3$ including the triples $(r_1,r_1\ldiv q,q)$
and $(q\rdiv c_1,c_1,q)$. Suppose the rows and columns of $\Sigma_3$ are
respectively $\{r_1,q\rdiv c_1,r_2\}$ and $\{r_1\ldiv q,c_1,c_2\}$.
Since $r_1,c_1\in S$ and $q\rdiv c_1\notin S$ we can see from
\lref{l:basic} that $\Sigma_0\cap\Sigma_3$ is a subsquare of order 1.
This means that $\{q\rdiv c_1,r_2,r_1\ldiv q,c_2\}\cap S=\emptyset$.
Also $r_1\circ(r_1\ldiv q)=(q\rdiv c_1)\circ c_1=r_2\circ c_2=q$.
It follows that $r_2\notin U_q$ and $c_2\notin T_q$.

Consider fixing $r_1$ and allowing $c_1$ to vary over $S$.  We find
that $r_2\circ(r_1\ldiv q)=r_1\circ c_1$ varies over $S$.  Next
allowing $r_1$ to vary over $S$, we conclude that the set of cells
$(Q\setminus\big(S\cup U_q)\big)\times T_q$ form a subsquare with symbols
$S$. Similarly, the cells $U_q\times\big(Q\setminus(S\cup T_q)\big)$ form a
subsquare with symbols $S$. The position of these subsquares is a
property of the loop independent of the choice of $q$. Hence for any
$q'\in Q\setminus S$ we must have either $U_{q'}=U_q$ and $T_{q'}=T_q$
or else $U_{q'}=Q\setminus(S\cup U_q)$ and $T_{q'}=Q\setminus(S\cup T_q)$.
It follows that the cells $U_q\times S$ and $S\times T_q$
form subsquares on the same symbols. Together with the subsquares
already identified, these are sufficient to show that $S$ is normal
in $Q$.
\end{proof}

\tref{t:normal} cannot be generalised to subloops of index $9$. For example,
in the non-abelian group of exponent $3$ and order $27$,
only $1$ of the $13$ subgroups of order $3$ is normal.

\begin{theorem}\label{t:min}
Let the loop $(Q,\cdot)$ be a minimal counterexample to \cjref{c:vR}.
Then $Q$ is simple.
\end{theorem}

\begin{proof}
If $Q$ had a proper normal subloop $N$, then by Corollary \ref{c:sublpclosed},
both $N$ and $Q/N$ would be van Rees loops. By minimality of $Q$, each of
$|N|$ and $|Q/N|$ are powers of $3$, and hence, so is $|Q|$, a contradiction.
\end{proof}

\section{Examples and Classification}
\label{s:constructions}

In this section we consider various examples of van Rees loops.
Obvious examples are elementary abelian $3$-groups and nonabelian
groups of exponent $3$. These can never yield a counterexample to
Conjecture~\ref{c:vR}.  Thus we are primarily interested in finding
nonassociative examples. We can obviously rule out order $3$, and, by
the following result, order $9$.

\begin{theorem}
A van Rees loop of order $9$ is an elementary abelian $3$-group.
\end{theorem}

\begin{proof}
Let $Q$ be a van Rees loop of order $9$.  Any nonidentity
element generates a subloop of order $3$, and every such
subloop is normal by \tref{t:normal}.  Take any two distinct such
subloops, say $H$ and $K$. Noting that $|H|\cdot |K| = 9$ and
$H\cap K = \{1\}$, we have that $Q$ is a direct product of cyclic
groups of order $3$. Thus $Q$ is associative and elementary abelian as
claimed.
\end{proof}

Recall that the \emph{left}, \emph{middle} and \emph{right nuclei} of a
loop $Q$ are the sets
\begin{align*}
N_{\lambda}(Q) &= \{a\in Q\ |\ (ax)y=a(xy),\ \forall x,y\in Q\}, \\
N_{\mu}(Q) &= \{a\in Q\ |\ (xa)y=x(ay),\ \forall x,y\in Q\},\\
N_{\rho}(Q) &= \{a\in Q\ |\ (xy)a=x(ya),\ \forall x,y\in Q\},
\end{align*}
respectively. The \emph{center}
$Z(Q) = N_{\lambda}(Q)\cap N_{\mu}(Q) \cap N_{\rho}(Q)
\cap \{a\in Q\ |\ ax=xa,\ \forall x\in Q\}$.

\setlength{\tabcolsep}{3pt}

\begin{table}[htb]
\begin{center}
\small\texttt{
\begin{tabular}{c|ccccccccccccccccccccccccccc}
$\cdot$&\id&a&b&c&d&e&f&g&h&i&j&k&l&m&n&o&p&q&r&s&t&u&v&w&x&y&z\\
\hline
\id&\id&a&b&c&d&e&f&g&h&i&j&k&l&m&n&o&p&q&r&s&t&u&v&w&x&y&z\\
b&b&\id&a&g&e&h&c&f&d&m&i&p&n&j&o&l&q&k&v&r&y&w&s&x&u&z&t\\
c&c&f&g&d&\id&b&h&e&a&l&o&i&k&n&p&q&m&j&u&x&r&t&w&y&z&v&s\\
d&d&h&e&\id&c&g&a&b&f&k&q&l&i&p&m&j&n&o&t&z&u&r&y&v&s&w&x\\
e&e&d&h&b&g&f&\id&a&c&p&k&n&m&q&j&i&o&l&y&t&w&v&z&s&r&x&u\\
f&f&g&c&h&a&\id&e&d&b&o&n&j&q&l&k&p&i&m&x&w&s&z&u&t&y&r&v\\
g&g&c&f&e&b&a&d&h&\id&n&l&m&p&o&q&k&j&i&w&u&v&y&x&z&t&s&r\\
h&h&e&d&a&f&c&b&\id&g&q&p&o&j&k&i&m&l&n&z&y&x&s&t&r&v&u&w\\
i&i&k&l&j&m&n&q&o&p&r&x&z&w&y&s&t&u&v&\id&e&g&h&f&b&c&d&a\\
j&j&q&o&m&i&l&p&n&k&t&s&x&v&w&z&u&r&y&f&\id&d&b&e&c&h&a&g\\
k&k&l&i&q&p&m&o&j&n&v&u&t&x&z&r&y&w&s&h&c&\id&g&b&d&a&f&e\\
l&l&i&k&o&n&p&j&q&m&s&w&y&u&t&v&z&x&r&g&a&h&\id&d&f&e&b&c\\
m&m&p&n&i&j&o&k&l&q&u&z&s&y&v&x&r&t&w&e&f&a&c&\id&h&b&g&d\\
n&n&m&p&l&o&q&i&k&j&x&y&v&t&r&w&s&z&u&a&d&c&e&g&\id&f&h&b\\
o&o&j&q&n&l&k&m&p&i&z&v&w&r&u&y&x&s&t&d&g&b&f&a&e&\id&c&h\\
p&p&n&m&k&q&j&l&i&o&w&t&r&z&s&u&v&y&x&c&b&e&a&h&g&d&\id&f\\
q&q&o&j&p&k&i&n&m&l&y&r&u&s&x&t&w&v&z&b&h&f&d&c&a&g&e&\id\\
r&r&x&y&w&z&u&t&s&v&\id&d&b&a&c&f&h&g&e&i&p&n&q&o&m&l&k&j\\
s&s&w&t&u&y&x&z&v&r&c&\id&g&f&d&h&a&e&b&n&j&q&i&k&o&p&m&l\\
t&t&s&w&v&x&r&u&z&y&a&h&\id&b&f&g&e&c&d&o&i&k&m&p&l&q&j&n\\
u&u&z&v&y&s&t&r&x&w&b&e&a&\id&g&c&d&f&h&p&o&j&l&i&q&m&n&k\\
v&v&u&z&x&t&w&y&r&s&d&c&e&h&\id&a&f&b&g&q&l&i&n&m&k&j&o&p\\
w&w&t&s&z&r&y&v&u&x&e&g&h&d&b&\id&c&a&f&l&m&p&j&q&n&k&i&o\\
x&x&y&r&t&v&z&s&w&u&g&b&f&c&e&d&\id&h&a&j&k&l&p&n&i&o&q&m\\
y&y&r&x&s&u&v&w&t&z&h&f&d&e&a&b&g&\id&c&m&q&o&k&l&j&n&p&i\\
z&z&v&u&r&w&s&x&y&t&f&a&c&g&h&e&b&d&\id&k&n&m&o&j&p&i&l&q\\
\end{tabular}}
\end{center}
\caption{A van Rees loop with all nuclei trivial.}
\label{t:trivnuc}
\end{table}

Loops $Q_1, Q_2$ are said to be \emph{paratopic} (or
\emph{isostrophic}) if $Q_1$ is isotopic to a conjugate (or
parastrophe) of the other.
Using a computer search, we have classified up to paratopy all van
Rees loops of order $27$ such that at least one of the nuclei is
nontrivial. (Note that it is irrelevant which nucleus is specified to
be nontrivial. If a loop has, say, nontrivial left nucleus, then there
is a paratope with nontrivial
middle nucleus and a paratope with nontrivial right nucleus.)  We should
stress that there are also examples of van Rees loops with all nuclei
trivial. We do not know how many, but \Tref{t:trivnuc} gives one example.

\begin{table}[htb]
\begin{center}
\small\texttt{
\begin{tabular}{c|ccccccccccccccccccccccccccc}
$\cdot$&\id&a&b&c&d&e&f&g&h&i&j&k&l&m&n&o&p&q&r&s&t&u&v&w&x&y&z\\
\hline
\id&\id&a&b&c&d&e&f&g&h&i&j&k&l&m&n&o&p&q&r&s&t&u&v&w&x&y&z\\
a&a&b&\id&e&g&f&c&h&d&j&k&i&m&n&l&p&q&o&s&t&r&v&w&u&y&z&x\\
b&b&\id&a&f&h&c&e&d&g&k&i&j&n&l&m&q&o&p&t&r&s&w&u&v&z&x&y\\
c&c&e&f&d&\id&g&h&a&b&o&p&q&r&s&t&w&u&v&y&z&x&j&k&i&n&l&m\\
d&d&g&h&\id&c&a&b&e&f&w&u&v&y&z&x&i&j&k&l&m&n&p&q&o&t&r&s\\
e&e&f&c&g&a&h&d&b&\id&p&q&o&s&t&r&u&v&w&z&x&y&k&i&j&l&m&n\\
f&f&c&e&h&b&d&g&\id&a&q&o&p&t&r&s&v&w&u&x&y&z&i&j&k&m&n&l\\
g&g&h&d&a&e&b&\id&f&c&u&v&w&z&x&y&j&k&i&m&n&l&q&o&p&r&s&t\\
h&h&d&g&b&f&\id&a&c&e&v&w&u&x&y&z&k&i&j&n&l&m&o&p&q&s&t&r\\
i&i&j&k&o&w&p&q&u&v&l&m&n&\id&a&b&r&s&t&c&e&f&z&x&y&h&d&g\\
j&j&k&i&p&u&q&o&v&w&m&n&l&a&b&\id&s&t&r&e&f&c&x&y&z&d&g&h\\
k&k&i&j&q&v&o&p&w&u&n&l&m&b&\id&a&t&r&s&f&c&e&y&z&x&g&h&d\\
l&l&m&n&r&y&s&t&z&x&\id&a&b&i&j&k&c&e&f&o&p&q&g&h&d&v&w&u\\
m&m&n&l&s&z&t&r&x&y&a&b&\id&j&k&i&e&f&c&p&q&o&h&d&g&w&u&v\\
n&n&l&m&t&x&r&s&y&z&b&\id&a&k&i&j&f&c&e&q&o&p&d&g&h&u&v&w\\
o&o&p&q&u&k&v&w&i&j&s&t&r&f&c&e&x&y&z&d&g&h&m&n&l&\id&a&b\\
p&p&q&o&v&i&w&u&j&k&t&r&s&c&e&f&y&z&x&g&h&d&n&l&m&a&b&\id\\
q&q&o&p&w&j&u&v&k&i&r&s&t&e&f&c&z&x&y&h&d&g&l&m&n&b&\id&a\\
r&r&s&t&z&n&x&y&l&m&f&c&e&p&q&o&d&g&h&v&w&u&b&\id&a&k&i&j\\
s&s&t&r&x&l&y&z&m&n&c&e&f&q&o&p&g&h&d&w&u&v&\id&a&b&i&j&k\\
t&t&r&s&y&m&z&x&n&l&e&f&c&o&p&q&h&d&g&u&v&w&a&b&\id&j&k&i\\
u&u&v&w&i&q&j&k&o&p&x&y&z&d&g&h&m&n&l&b&\id&a&s&t&r&c&e&f\\
v&v&w&u&j&o&k&i&p&q&y&z&x&g&h&d&n&l&m&\id&a&b&t&r&s&e&f&c\\
w&w&u&v&k&p&i&j&q&o&z&x&y&h&d&g&l&m&n&a&b&\id&r&s&t&f&c&e\\
x&x&y&z&m&r&n&l&s&t&g&h&d&w&u&v&\id&a&b&k&i&j&c&e&f&o&p&q\\
y&y&z&x&n&s&l&m&t&r&h&d&g&u&v&w&a&b&\id&i&j&k&e&f&c&p&q&o\\
z&z&x&y&l&t&m&n&r&s&d&g&h&v&w&u&b&\id&a&j&k&i&f&c&e&q&o&p
\end{tabular}}
\end{center}
\caption{A universal left conjugacy closed, left inverse property, van Rees loop.}
\label{t:lcc}
\end{table}

\begin{theorem}\label{t:order27}
Up to paratopy, there are exactly six nonassociative van Rees
loops of order $27$ with at least one nontrivial nucleus.
\end{theorem}

The six species in the theorem include the following representatives 
\cite{wwww}, each in a different class:
\begin{itemize}
\item A Bol loop with trivial center, discovered by Keedwell
\cite{Keed0,Keed1} and described in \cite{FK}.
\item Two power-associative conjugacy closed loops, described in \cite{KK}.
\item A universal left conjugacy closed loop (which is not conjugacy closed) with the
left inverse property; see Table \ref{t:lcc}.
\item A commutative, weak inverse property loop; see Table \ref{t:commwip27}.
\item A (noncommutative) weak inverse property loop such that each inner mapping
of the form $L_x^{-1} R_x$ is an automorphism; see Table \ref{t:noncommwip}.
\end{itemize}
The species with the Bol loop is the only one such that each loop in
the species has trivial center. It, and the species of the
universal left conjugacy closed loop both contain 3 distinct isotopy classes.
The other 4 species contain a single isotopy class, so there are 10 isotopy
classes of nonassociative van Rees loops of order 27 with a nontrivial nucleus.

The Bol loop in \tref{t:order27} and the loop in \Tref{t:trivnuc} each have
a single subloop of order 9. By \tref{t:normal} the subloops are normal,
leading to 9 latin subsquares of order 9 in each Cayley table.
The 5 loops in \tref{t:order27} other than the bol loop each
have 4 (normal) subloops of order 9, hence 36 subsquares of order 9 in
their Cayley table.

\begin{table}[htb]
\begin{center}
\small\texttt{
\begin{tabular}{c|ccccccccccccccccccccccccccc}
$\cdot$&\id&a&b&c&d&e&f&g&h&i&j&k&l&m&n&o&p&q&r&s&t&u&v&w&x&y&z\\
\hline
\id&\id&a&b&c&d&e&f&g&h&i&j&k&l&m&n&o&p&q&r&s&t&u&v&w&x&y&z\\
a&a&b&\id&e&g&f&c&h&d&j&k&i&m&n&l&p&q&o&s&t&r&v&w&u&y&z&x\\
b&b&\id&a&f&h&c&e&d&g&k&i&j&n&l&m&q&o&p&t&r&s&w&u&v&z&x&y\\
c&c&e&f&d&\id&g&h&a&b&o&p&q&t&r&s&u&v&w&y&z&x&i&j&k&l&m&n\\
d&d&g&h&\id&c&a&b&e&f&w&u&v&z&x&y&k&i&j&l&m&n&q&o&p&s&t&r\\
e&e&f&c&g&a&h&d&b&\id&p&q&o&r&s&t&v&w&u&z&x&y&j&k&i&m&n&l\\
f&f&c&e&h&b&d&g&\id&a&q&o&p&s&t&r&w&u&v&x&y&z&k&i&j&n&l&m\\
g&g&h&d&a&e&b&\id&f&c&u&v&w&x&y&z&i&j&k&m&n&l&o&p&q&t&r&s\\
h&h&d&g&b&f&\id&a&c&e&v&w&u&y&z&x&j&k&i&n&l&m&p&q&o&r&s&t\\
i&i&j&k&o&w&p&q&u&v&l&m&n&\id&a&b&r&s&t&c&e&f&y&z&x&d&g&h\\
j&j&k&i&p&u&q&o&v&w&m&n&l&a&b&\id&s&t&r&e&f&c&z&x&y&g&h&d\\
k&k&i&j&q&v&o&p&w&u&n&l&m&b&\id&a&t&r&s&f&c&e&x&y&z&h&d&g\\
l&l&m&n&t&z&r&s&x&y&\id&a&b&i&j&k&f&c&e&q&o&p&d&g&h&v&w&u\\
m&m&n&l&r&x&s&t&y&z&a&b&\id&j&k&i&c&e&f&o&p&q&g&h&d&w&u&v\\
n&n&l&m&s&y&t&r&z&x&b&\id&a&k&i&j&e&f&c&p&q&o&h&d&g&u&v&w\\
o&o&p&q&u&k&v&w&i&j&r&s&t&f&c&e&x&y&z&g&h&d&m&n&l&\id&a&b\\
p&p&q&o&v&i&w&u&j&k&s&t&r&c&e&f&y&z&x&h&d&g&n&l&m&a&b&\id\\
q&q&o&p&w&j&u&v&k&i&t&r&s&e&f&c&z&x&y&d&g&h&l&m&n&b&\id&a\\
r&r&s&t&y&l&z&x&m&n&c&e&f&q&o&p&g&h&d&w&u&v&a&b&\id&k&i&j\\
s&s&t&r&z&m&x&y&n&l&e&f&c&o&p&q&h&d&g&u&v&w&b&\id&a&i&j&k\\
t&t&r&s&x&n&y&z&l&m&f&c&e&p&q&o&d&g&h&v&w&u&\id&a&b&j&k&i\\
u&u&v&w&i&q&j&k&o&p&y&z&x&d&g&h&m&n&l&a&b&\id&t&r&s&f&c&e\\
v&v&w&u&j&o&k&i&p&q&z&x&y&g&h&d&n&l&m&b&\id&a&r&s&t&c&e&f\\
w&w&u&v&k&p&i&j&q&o&x&y&z&h&d&g&l&m&n&\id&a&b&s&t&r&e&f&c\\
x&x&y&z&l&s&m&n&t&r&d&g&h&v&w&u&\id&a&b&k&i&j&f&c&e&o&p&q\\
y&y&z&x&m&t&n&l&r&s&g&h&d&w&u&v&a&b&\id&i&j&k&c&e&f&p&q&o\\
z&z&x&y&n&r&l&m&s&t&h&d&g&u&v&w&b&\id&a&j&k&i&e&f&c&q&o&p
\end{tabular}}
\end{center}
\caption{A commutative, weak inverse property, van Rees loop.}
\label{t:commwip27}
\end{table}

\begin{table}[htb]
\begin{center}
\small\texttt{
\begin{tabular}{c|ccccccccccccccccccccccccccc}
$\cdot$&\id&a&b&c&d&e&f&g&h&i&j&k&l&m&n&o&p&q&r&s&t&u&v&w&x&y&z\\
\hline
\id&\id&a&b&c&d&e&f&g&h&i&j&k&l&m&n&o&p&q&r&s&t&u&v&w&x&y&z\\
a&a&b&\id&e&g&f&c&h&d&j&k&i&m&n&l&p&q&o&s&t&r&v&w&u&y&z&x\\
b&b&\id&a&f&h&c&e&d&g&k&i&j&n&l&m&q&o&p&t&r&s&w&u&v&z&x&y\\
c&c&e&f&d&\id&g&h&a&b&o&p&q&s&t&r&u&v&w&y&z&x&i&j&k&m&n&l\\
d&d&g&h&\id&c&a&b&e&f&v&w&u&x&y&z&j&k&i&l&m&n&p&q&o&r&s&t\\
e&e&f&c&g&a&h&d&b&\id&p&q&o&t&r&s&v&w&u&z&x&y&j&k&i&n&l&m\\
f&f&c&e&h&b&d&g&\id&a&q&o&p&r&s&t&w&u&v&x&y&z&k&i&j&l&m&n\\
g&g&h&d&a&e&b&\id&f&c&w&u&v&y&z&x&k&i&j&m&n&l&q&o&p&s&t&r\\
h&h&d&g&b&f&\id&a&c&e&u&v&w&z&x&y&i&j&k&n&l&m&o&p&q&t&r&s\\
i&i&j&k&p&u&q&o&v&w&l&m&n&\id&a&b&s&t&r&e&f&c&x&y&z&d&g&h\\
j&j&k&i&q&v&o&p&w&u&m&n&l&a&b&\id&t&r&s&f&c&e&y&z&x&g&h&d\\
k&k&i&j&o&w&p&q&u&v&n&l&m&b&\id&a&r&s&t&c&e&f&z&x&y&h&d&g\\
l&l&m&n&r&y&s&t&z&x&\id&a&b&i&j&k&c&e&f&o&p&q&g&h&d&v&w&u\\
m&m&n&l&s&z&t&r&x&y&a&b&\id&j&k&i&e&f&c&p&q&o&h&d&g&w&u&v\\
n&n&l&m&t&x&r&s&y&z&b&\id&a&k&i&j&f&c&e&q&o&p&d&g&h&u&v&w\\
o&o&p&q&v&i&w&u&j&k&r&s&t&e&f&c&x&y&z&d&g&h&m&n&l&\id&a&b\\
p&p&q&o&w&j&u&v&k&i&s&t&r&f&c&e&y&z&x&g&h&d&n&l&m&a&b&\id\\
q&q&o&p&u&k&v&w&i&j&t&r&s&c&e&f&z&x&y&h&d&g&l&m&n&b&\id&a\\
r&r&s&t&x&m&y&z&n&l&c&e&f&p&q&o&g&h&d&u&v&w&\id&a&b&i&j&k\\
s&s&t&r&y&n&z&x&l&m&e&f&c&q&o&p&h&d&g&v&w&u&a&b&\id&j&k&i\\
t&t&r&s&z&l&x&y&m&n&f&c&e&o&p&q&d&g&h&w&u&v&b&\id&a&k&i&j\\
u&u&v&w&j&o&k&i&p&q&y&z&x&d&g&h&l&m&n&\id&a&b&r&s&t&e&f&c\\
v&v&w&u&k&p&i&j&q&o&z&x&y&g&h&d&m&n&l&a&b&\id&s&t&r&f&c&e\\
w&w&u&v&i&q&j&k&o&p&x&y&z&h&d&g&n&l&m&b&\id&a&t&r&s&c&e&f\\
x&x&y&z&l&s&m&n&t&r&g&h&d&u&v&w&\id&a&b&j&k&i&c&e&f&o&p&q\\
y&y&z&x&m&t&n&l&r&s&h&d&g&v&w&u&a&b&\id&k&i&j&e&f&c&p&q&o\\
z&z&x&y&n&r&l&m&s&t&d&g&h&w&u&v&b&\id&a&i&j&k&f&c&e&q&o&p
\end{tabular}}
\end{center}
\caption{A noncommutative, weak inverse property, van Rees loop with each $L_x^{-1} R_x$ an automorphism.}
\label{t:noncommwip}
\end{table}

Finally, we mention that in certain varieties of loops, any loop 
of exponent $3$ is necessarily a van Rees loop.

\smallskip

\noindent\emph{Conjugacy closed (CC-)loops}: Every CC-loop is isomorphic
to all of its isotopes \cite{GR}. Thus a CC-loop of exponent $3$
is a van Rees loop.

\smallskip

\noindent\emph{Moufang loops}: It is classical that every Moufang loop of exponent
$3$ is a van Rees loop \cite{Bruck}. More generally, this is subsumed by the following.

\smallskip

\noindent\emph{Bol loops}: For a Bol loop of odd exponent, the exponent is an isotopy
invariant (\cite{Hall}, Corollary 6.7). Thus a Bol loop of exponent $3$ is a van Rees
loop.

\smallskip

The three varieties of loops just discussed are all \emph{isotopically invariant}, that is,
any loop isotope of a loop in the variety is also in the variety.  A natural question arises:

\smallskip

\noindent\emph{In which isotopically invariant varieties of loops must a loop of
exponent} $3$ \emph{necessarily be a van Rees loop?}

\smallskip

\noindent The general answer here, of course, is ``not all varieties''. 
We have already seen in \egref{ex:notinv} that exponent $3$
is not intrinsically an isotopically invariant property.
So the answer to the question is no for the variety of \emph{all} loops.
Similarly, while power-associativity is not an isotopically invariant property, one might
consider the variety of loops such that every isotope is power-associative. 
Example~\ref{ex:notinv} once again shows that the answer is no in this case.

\smallskip

We conclude this section with another open problem motivated by our
examples.  A loop is \emph{diassociative} if every $2$-generated
subloop is a group. For example, a consequence of Moufang's Theorem is
that Moufang loops are diassociative, but non-Moufang diassociative
loops exist as well. Interestingly, the only diassociative van Rees
loops known to us are Moufang loops of exponent $3$.

\begin{problem}
Let $Q$ be a diassociative van Rees loop. Must $Q$ be a Moufang loop?
\end{problem}

\section{Steiner quasigroups}\label{s:steiner}

Recall that a quasigroup $(Q,\cdot)$ is said to be \emph{Steiner}
if it satisfies the identities $xx = x$, $x\cdot yx = y$ and $xy = yx$, 
that is, a Steiner quasigroup is idempotent and totally symmetric.
Steiner quasigroups are essentially the same as \emph{Steiner triple
systems} in that the triples of the latter define the multiplication
of the former.

If $(Q,\cdot)$ is a Steiner quasigroup, then fixing an element $\id\in Q$,
we can define a loop isotope by $x\circ y = (\id x)(\id y)$. Then
$(Q,\circ)$ is commutative and has exponent $3$, and also satisfies the
weak inverse property (WIP) $x\circ ((y\circ x)\circ (y\circ x)) = y\circ y$.
Set $\mathcal{W}(Q,\cdot) = (Q,\circ)$. On the other hand,
if we start with a commutative WIP loop $(Q,\circ)$ of exponent $3$ and
define $x\ast y = (x\circ x)\circ (y\circ y)$, we obtain a Steiner 
quasigroup $(Q,\ast)$.
Set $\mathcal{S}(Q,\circ) = (Q,\ast)$. It is easy to check that 
$\mathcal{S}\mathcal{W}(Q,\cdot) = (Q,\cdot)$ and 
$\mathcal{W}\mathcal{S}(Q,\circ) = (Q,\circ)$, provided that the fixed
element of $Q$ we use to construct $\mathcal{W}\mathcal{S}(Q,\circ)$
is the identity element of $(Q,\circ)$. Summing up, the varieties of
Steiner quasigroups and commutative WIP loops of exponent $3$ are 
term equivalent. 

Referring to Theorem~\ref{t:variety}, we see that a Steiner quasigroup is a
van Rees quasigroup if and only if it satisfies the identity 
\[
x(y\cdot xz) = y(x\cdot yz).
\]
For loops, the following is an interesting aside.

\begin{theorem}
A commutative van Rees loop has the weak inverse property.
\end{theorem}

\begin{proof}
Firstly, we compute 
\[
x^2 (y\cdot yx) = x^2 (y\cdot y(x^2 x^2)) = y(x^2\cdot x^2 (y x^2)) 
= y(x^2\cdot x^2 (x^2 y)) = y^2,
\]
using (vRL3) and (vRL1). Now replacing $x$ with $yx$, we have 
\[
y^2 = (yx)^2 (y\cdot y(yx)) = (yx)^2 x = x(yx)^2,
\]
using (vRL1) once more, and commutativity. This completes the proof.
\end{proof}

\emph{Distributive} Steiner quasigroups are defined by the identity 
$x\cdot yz = xy\cdot xz$. The corresponding loop isotopes as described
above are commutative Moufang loops of exponent $3$. Thus distributive
Steiner quasigroups are certainly van Rees quasigroups, although this
can just as easily be seen directly: 
$x(y\cdot xz) = xy\cdot (x\cdot
xz) = xy\cdot z = yx\cdot z = yx\cdot (y\cdot yz) = y(x\cdot yz)$.

To see that van Rees Steiner quasigroups arise naturally among nondistributive
quasigroups, we note Marczak's investigation of essentially ternary polynomials
in nondistributive Steiner quasigroups \cite{Mar97}. He showed that for a
such a quasigroup, there are at least $21$ essentially ternary polynomials,
and that this bound is achieved if and only if the following identity 
holds:
\[
(xz\cdot yz)\cdot (xy\cdot z) = x(y\cdot xz).  \tag{M}
\]
We will call Steiner quasigroups satisfying (M) \emph{Marczak} Steiner
quasigroups. It is easy to see that distributive Steiner quasigroups are
Marczak. On the other hand, the left hand side of (M) is clearly invariant
under exchanging the variables $x$ and $y$, and so the same is true of the
right side. Thus every Marczak Steiner quasigroup is a van Rees quasigroup.
We thus have the following inclusions of varieties of Steiner quasigroups:
\[
\text{Distributive}\qquad \subset \qquad \text{Marczak}\qquad \subset \qquad \text{van Rees}.
\]
The first inclusion is proper; indeed, Marczak showed that (in our terminology)
any nondistributive, Marczak Steiner quasigroup must contain a certain subquasigroup
of order $27$. That quasigroup is among our examples; its commutative WIP loop
isotope is the loop given by Table~\ref{t:commwip27}. 

The second inclusion above is also proper. We omit the table (although
it is available at \cite{wwww}), but
using \textsc{Mace4} \cite{McCune}, we have found a non-Marczak, van
Rees Steiner quasigroup of order $81$.  It is a bit easier to describe
its loop isotope $Q$, which is nilpotent of class $3$.  The center
$Z(Q)$ has order $3$ and the factor $Q/Z(Q)$ is isomorphic to the loop
in Table~\ref{t:commwip27}. The second center $Z_2(Q)$ has order $9$
and coincides with the derived subloop.

We do not know if van Rees' Conjecture holds for van Rees Steiner
quasigroups.  An affirmative answer to \cjref{c:vR} in the
Steiner case would generalize the well-known fact that every finite
distributive Steiner quasigroup has order a power of $3$.

\begin{problem}
If there exists a Steiner counterexample to \cjref{c:vR}, is it the case that
every finite Marczak Steiner quasigroup has order a power of $3$?
\end{problem}

\medskip

In the remainder of this section, we describe a construction associating to any
van Rees quasigroup or loop a family of Steiner quasigroups
on the same underlying set. 

Suppose that $L$ is a van Rees latin square of order $n$.
Each row, column or symbol of $L$ induces a Steiner triple system of
order $n$.  We describe the construction for a given symbol $x$, which
we call the \emph{key}, but a
similar idea works for rows or columns.  The points of our Steiner
triple system are the occurrences of $x$ within $L$, and the triples
are the sets of $3$ points that are contained within a subsquare of
order $3$.  By \tref{t:main}(\ref{i:eachpair}) each pair of points
lies in a unique triple, as required. Applying this construction,
(regardless of whether the key is a row, column, or symbol), to a
van Rees loop $Q$ yields a Steiner triple system, which is equivalent
to a Steiner quasigroup in the usual way.

We can also describe the construction more algebraically, this time
using a given column $c$ as the key. The process is illustrated
in \eref{e:colSTS} below.
Two general symbols $a$ and $b$ occurring in column $c$ 
will occur in rows $a\rdiv c$ and $b\rdiv c$, respectively. 
In row $a\rdiv c$, $b$ occurs in column $(a\rdiv c)\ldiv b$.
In row $b\rdiv c$, $a$ occurs in column $(b\rdiv c)\ldiv a$. 
Thus the third element of the triple containing $a$ and $b$ is
$t=(a\rdiv c)[(b\rdiv c)\ldiv a] = (b\rdiv c)[(a\rdiv c)\ldiv b]$. 
\begin{equation}\label{e:colSTS}
\begin{array}{c|ccc}
& \ c \ & (a\rdiv c)\ldiv b & (b\rdiv c)\ldiv a \\
\hline
a\rdiv c & a & b & t\\
b\rdiv c & b & t & a\\
\end{array}
\end{equation}
Summarizing, we have the following.

\begin{theorem}\label{t:steiner}
Let $Q$ be a van Rees loop. For $c\in Q$, define a binary
operation $\ast_c$ on $Q$ by
\[
x \ast_c y = (x\rdiv c)[(y\rdiv c)\ldiv x]
\]
for all $x,y\in Q$. Then $(Q,\ast_c)$ is a Steiner quasigroup.
\end{theorem}

For instance, for $c = \id$, we obtain $x\ast_{\id} y = x(y\ldiv x) =
x(y\cdot yx)$. Of course, using a row as the key yields a dual result,
with the operation just being the mirror image of that in
\tref{t:steiner}.

Let $(Q,\cdot)$ be the commutative WIP loop given in
\Tref{t:commwip27}.  Using any row, column or symbol as the key, if we
apply the construction to $(Q,\cdot)$, we obtain a quasigroup
isomorphic to $\mathcal{S}(Q,\cdot)$, which is isotopic to $(Q,\cdot)$.

For all other van Rees loops of order $27$ described in
\sref{s:constructions} and for all choices of the key, the
construction in \tref{t:steiner} produces a quasigroup corresponding
to the affine Steiner triple system $AG(3,3)$. Notably then, for each
loop, the Steiner quasigroups produced by different choices of the key
all turn out to be isomorphic to each other. We do not know if that is
a general phenomenon.

\begin{problem}
Let $Q$ be a van Rees loop. For distinct $c_1, c_2\in Q$, are the Steiner
quasigroups $(Q,\ast_{c_1})$ and $(Q,\ast_{c_2})$ isomorphic? If not, are
they isotopic?
\end{problem}

Further, it follows that in all cases we have seen, the Steiner quasigroups
associated to a van Rees loop are themselves van Rees quasigroups. Nevertheless,
we have not been able to resolve the following.

\begin{problem}\label{p:steiner}
Let $(Q,\cdot)$ be a van Rees loop and fix $c\in Q$. Is
$(Q,\ast_{c})$ a van Rees quasigroup?
\end{problem}

In fact, \cjref{c:vR} can be resolved affirmatively in general
if it is resolved in the Steiner case, and if \pref{p:steiner}
has an affirmative answer.

\medskip

Finally, we comment on another feature of our construction.  For the
Steiner triple system using the symbol $x$ as key, each symbol $y\ne
x$ induces a parallel class of the triple system.  This is because
there are $n/3$ disjoint subsquares of order $3$ that contain both
symbol $x$ and symbol $y$.  Applying this idea to all choices of $y$
gives a set of parallel classes that includes every triple exactly
twice. In design theory terms, this might be called a 
\emph{double resolution}. We are not aware of such a concept
yet being introduced into the literature, but it is just possible
that it may prove helpful in resolving \cjref{c:vR} or \pref{p:steiner}.

\end{document}